\documentclass[12pt]{article}
\usepackage[psamsfonts]{amssymb}
\usepackage{amsthm,amsmath}
\usepackage{graphicx}

\title{Counting non-isomorphic maximal independent sets of the $n$-cycle graph}

\author{%
Raymond Bisdorff \\
\small{Computer Science and Communications Research Unit, University of Luxembourg}\\
\small{162A, avenue de la Fa\"{\i}encerie, L-1511 Luxembourg, Luxembourg}\\
\small{raymond.bisdorff[at]uni.lu}%
\and %
Jean-Luc Marichal\footnote{Corresponding author.} \\
\small{Mathematics Research Unit, University of Luxembourg}\\
\small{162A, avenue de la Fa\"{\i}encerie, L-1511 Luxembourg, Luxembourg}\\
\small{jean-luc.marichal[at]uni.lu}
}

\date{\small Revised, November 25, 2008}

\begin{document}
\maketitle

\theoremstyle{plain}
\newtheorem{theorem}{Theorem}[section]
\newtheorem{lemma}{Lemma}[section]
\newtheorem{proposition}{Proposition}[section]
\newtheorem{corollary}{Corollary}[section]

\theoremstyle{definition}
\newtheorem{definition}{Definition}[section]
\newtheorem{example}{Example}[section]

\theoremstyle{remark}
\newtheorem*{conjecture}{\indent Conjecture}
\newtheorem*{remark}{\indent Remark}

\newcommand{\N}{\mathbb{N}}                     
\newcommand{\Z}{\mathbb{Z}}                     
\newcommand{\R}{\mathbb{R}}                     
\newcommand{\Vspace}{\vspace{2ex}}              
\newcommand{\Orb}{\mathrm{Orb}}
\newcommand{\orb}{\mathrm{orb}}
\newcommand{\Stab}{\mathrm{Stab}}
\newcommand{\Ocal}{\mathcal{O}}
\newcommand{\Xbf}{\mathbf{X}}

\begin{abstract}
The number of maximal independent sets of the $n$-cycle graph $C_n$ is known to be the $n$th term of the Perrin sequence. The action of the
automorphism group of $C_n$ on the family of these maximal independent sets partitions this family into disjoint orbits, which represent the
non-isomorphic (i.e., defined up to a rotation and a reflection) maximal independent sets. We provide exact formulas for the total number of
orbits and the number of orbits having a given number of isomorphic representatives. We also provide exact formulas for the total number of
unlabeled (i.e., defined up to a rotation) maximal independent sets and the number of unlabeled maximal independent sets having a given number
of isomorphic representatives. It turns out that these formulas involve both Perrin and Padovan sequences.
\end{abstract}

\medskip


\medskip


\section{Introduction}

Let $G=(V,E)$ be a simple undirected graph, with vertex set $V$ and edge set $E$. Recall that an \emph{independent set}\/ of $G$ is a subset $X$
of $V$ such that no two vertices in $X$ are adjacent. An independent set is \emph{maximal}\/ if it is not a proper subset of any other
independent set.

In this paper we are concerned with maximal independent sets (MISs) of the $n$-cycle $C_n$, i.e., the graph consisting of a cycle with $n$
vertices. F\"uredi \cite{Fur87} observed that the total number of MISs of $C_n$ is given by the $n$th term $p(n)$ of the Perrin sequence
(Sloane's A001608 \cite{Slo}), which is defined recursively as $p(1)=0$, $p(2)=2$, $p(3)=3$, and
$$
p(n)=p(n-2)+p(n-3) \qquad (n\geqslant 4).
$$

The action of the automorphism group $\mathrm{Aut}(C_n)$ of $C_n$ on the family of MISs of $C_n$ gives rise to a partition of this family into
orbits, each containing \emph{isomorphic}\/ (i.e., identical up to a rotation and a reflection) MISs. For instance, it is easy to verify that
$C_6$ has $p(6)=5$ MISs, which can be grouped into two orbits: an orbit with two isomorphic MISs of size 3 and an orbit with three isomorphic
MISs of size 2. Thus, $C_6$ has essentially two non-isomorphic MISs.

We give here a closed form formula for the number of orbits, that is, the number of non-isomorphic MISs of $C_n$ (see Theorem~\ref{thm:main}).
We also provide a formula for the number of orbits having a given number of isomorphic representatives. Finally, by modifying the concept of
isomorphic MISs, requiring that they be identical up to a rotation only (which means that their \emph{unlabeled}\/ versions are identical), we
provide a formula for the number of these unlabeled MISs (see Theorem~\ref{thm:mainSigma}) and the number of unlabeled MISs having a given
number of isomorphic representatives. Incidentally, these formulas involve both Perrin and Padovan sequences, which are constructed from the
same recurrence equation.

Although many results were obtained on the maximum number of MISs in a general graph (see for instance Chang and Jou \cite{ChaJou99} and Ying et
al.\ \cite{YinMenSagVat06} for recent results) it seems that the searching for the exact number of distinct MISs, or non-isomorphic MISs, in
special subclasses of graphs has received little attention (see however Euler \cite{Eul05} and Kitaev \cite{Kit06}). We hope that by solving the
case of cycle graphs we might bring more interest to this kind of challenging enumerative combinatorics problems.

\section{Some properties of $\mathrm{Aut}(C_n)$}\label{sec:autCn}

In this section we recall and give some properties of the automorphism group $\mathrm{Aut}(C_n)$ that we will use in the rest of the paper. More
information on automorphism groups of graphs can be found for instance in Godsil and Royle~\cite{GodRoy01}.

Denote by $V_n=\{v_0,v_1,\ldots,v_{n-1}\}$ the set of vertices of $C_n$ labeled either clockwise or counterclockwise. The index set of the
vertices is assumed to be the cyclic group $\Z_n=\{0,1,\ldots,n-1\}$ with addition modulo $n$.

It is well known that the group $\mathrm{Aut}(C_n)$ is isomorphic to the \emph{dihedral group}\/ of order $2n$,
$$
D_{2n}=\{1,\sigma,\ldots,\sigma^{n-1},\tau,\sigma\tau,\ldots,\sigma^{n-1}\tau\},
$$
where $\sigma:2^{V_n}\to 2^{V_n}$ is the rotation which carries every $v_i$ into $v_{i+1}$ and $\tau:2^{V_n}\to 2^{V_n}$ is the reflection which
carries every $v_i$ into $v_{n-i}$, with the properties $\sigma^n=1$, $\tau^2=1$, and $\tau\sigma=\sigma^{-1}\tau$. The subset
$\{1,\sigma,\ldots,\sigma^{n-1}\}$ of rotations, also denoted $\langle\sigma\rangle$, is a cyclic subgroup of order $n$ of $D_{2n}$.

Consider a family $\Xbf_n\subseteq 2^{V_n}$ of subsets of $V_n$ with the property that if $X\in\Xbf_n$ then $g(X)\in\Xbf_n$ for all $g\in
D_{2n}$. For instance, the family of MISs of $C_n$ fulfills this property. Indeed, independence and maximality properties, which characterize
the MISs, remain stable under the action of $D_{2n}$.

Recall that, for any $X\in \Xbf_n$, the \emph{orbit}\/ and the \emph{stabilizer}\/ of $X$ under the action of $D_{2n}$ are respectively defined
as
\begin{eqnarray*}
\Orb(X)     &=& \{g(X):g\in D_{2n}\},\\
\Stab(X)    &=& \{g\in D_{2n}:g(X)=X\}.
\end{eqnarray*}
Recall also that $\Stab(X)$ is either a cyclic or a dihedral subgroup of $D_{2n}$ and that, by the orbit-stabilizer theorem, we have
\begin{equation}\label{eq:OrbStabTh}
|\Orb(X)|\times |\Stab(X)|=|D_{2n}|=2n,
\end{equation}
which implies that both $|\Orb(X)|$ and $|\Stab(X)|$ divide $2n$.

Let $\Ocal_n=\Xbf_n/D_{2n}$ denote the set of orbits of $\Xbf_n$ under the action of $D_{2n}$. We say that $X,X'\in \Xbf_n$ are
\emph{isomorphic} if $\Orb(X)=\Orb(X')$, that is, if $X'$ can be obtained from $X$ by a rotation and/or a reflection.

For any divisor $d\geqslant 1$ of $2n$ (we write $d|2n$), denote by $\Ocal_n^d$ the set of orbits of $\Xbf_n$ of cardinality $\frac{2n}d$. By
(\ref{eq:OrbStabTh}), $\Ocal_n^d$ is also the set of orbits whose elements have a stabilizer of cardinality $d$. Define also the sequences
\begin{eqnarray*}
\orb(n) &=& |\Ocal_n|,\\
\orb_d(n) &=& |\Ocal_n^d| \qquad (d|2n),
\end{eqnarray*}
which obviously lead to the identity
\begin{equation}\label{eq:NumberOrbits}
\orb(n) = \sum_{d|2n} \orb_d(n).
\end{equation}
Moreover, since the orbits partition the set $\Xbf_n$ we immediately have
\begin{equation}\label{eq:XDecomposition}
|\Xbf_n| = \sum_{d|2n} \frac{2n}d\, \orb_d(n).
\end{equation}

We now introduce a concept that will be very useful as we continue, namely the \emph{membership function}\/ of any set $X\in\Xbf_n$.

\begin{definition}
The \emph{membership}\/ (or \emph{characteristic}) \emph{function}\/ of a set $X\in\Xbf_n$ is the mapping $\mathbf{1}_X:V\to\{0,1\}$ defined as
$\mathbf{1}_X(v_i)=1$ if and only if $v_i\in X$.
\end{definition}

The following immediate result expresses, for any set $X\in\Xbf_n$, the properties of $\Stab(X)$ in terms of the membership function
$\mathbf{1}_X$.

\begin{proposition}\label{prop:RotRef}
For any $j\in\Z$ and any $X\in\Xbf_n$, we have
\begin{itemize}
\item[(i)] $\sigma^{j}\in\Stab(X)$ if and only if $\mathbf{1}_X(v_{\ell +j})=\mathbf{1}_X(v_{\ell})$ for all $\ell\in\Z$,

\item[(ii)] $\sigma^{2j}\tau\in\Stab(X)$ if and only if $\mathbf{1}_X(v_{j+\ell})=\mathbf{1}_X(v_{j-\ell})$ for all $\ell\in\Z$,

\item[(iii)] $\sigma^{2j+1}\tau\in\Stab(X)$ if and only if $\mathbf{1}_X(v_{j+\ell})=\mathbf{1}_X(v_{j+1-\ell})$ for all $\ell\in\Z$.
\end{itemize}
\end{proposition}

The conditions stated in Proposition~\ref{prop:RotRef} have a clear geometric interpretation. If we think of $C_n$ as a regular $n$-gon centered
at the origin, then condition $(i)$ means that $X$ is invariant under the rotation, about the origin, which carries $v_0$ into $v_j$. Then,
condition $(ii)$ means that $X$ has a symmetry axis passing through $v_j$. Finally, condition $(iii)$ means that $X$ has a symmetry axis passing
through the midpoint of $v_j$ and $v_{j+1}$.

From Proposition~\ref{prop:RotRef} it follows that, if $n$ is odd, all reflections fix exactly one vertex while, in the case of $n$ even, all
reflections of the form $\sigma^i\tau$ fix exactly two vertices or no vertex according to whether $i$ is even or odd, respectively.

\begin{definition}\label{de:dConc}
Given $X\in \Xbf_n$ and an integer $d\geqslant 1$, the \emph{$d$-concatenation}\/ of $X$ is the set $X^{(d)}\in\Xbf_{dn}$ whose membership
function is given by
$$
\mathbf{1}_{X^{(d)}}(v_{\ell})=\mathbf{1}_{X}(v_{\ell})\qquad(\ell\in \Z).
$$
\end{definition}

Definition~\ref{de:dConc} can be reformulated as follows: For any $X\in \Xbf_n$, the set $X^{(d)}\in\Xbf_{dn}$ is characterized by the fact that
the vector $\mathbf{1}_{X^{(d)}}(V_{dn})$ is obtained by concatenating $d$ times the vector $\mathbf{1}_X(V_n)$.

\begin{proposition}\label{prop:XXd}
For any $X\in \Xbf_n$ and any $j\in\Z$, we have $\sigma^j\in\Stab(X)$ if and only if $\sigma^j\in\Stab(X^{(d)})$. Similarly, we have
$\sigma^j\tau\in\Stab(X)$ if and only if $\sigma^j\tau\in\Stab(X^{(d)})$.
\end{proposition}

\begin{proof}
Let us prove the second part. The first part can be proved similarly. Let $X\in \Xbf_n$ and $j\in\Z$. Then
\begin{eqnarray*}
\sigma^j\tau\in\Stab(X)
& \Leftrightarrow & \mathbf{1}_X(v_{j-\ell})=\mathbf{1}_X(v_{\ell})\qquad (\ell\in\Z)\\
& \Leftrightarrow & \mathbf{1}_{X^{(d)}}(v_{j-\ell})=\mathbf{1}_{X^{(d)}}(v_{\ell})\qquad (\ell\in\Z)\\
& \Leftrightarrow & \sigma^j\tau\in\Stab(X^{(d)}).
\end{eqnarray*}
\end{proof}

\begin{proposition}\label{prop:NOrbIncl}
The number of orbits of $\Xbf_n$ whose elements have a stabilizer included in $\langle\sigma\rangle$ is given by $\sum_{d|n}\orb_1(d)$.
\end{proposition}

\begin{proof}
Let $d|n$ and let $X\in\Xbf_n$ be such that $\Stab(X)=\langle\sigma^d\rangle$. Then
$$
\mathbf{1}_X(v_{\ell+id})=\mathbf{1}_X(v_{\ell}) \qquad (i,\ell\in\Z)
$$
and hence, there is $Y\in\Xbf_d$ such that $X=Y^{(n/d)}$. By Proposition~\ref{prop:XXd}, we necessarily have $\Stab(Y)=\{1\}$ and hence
$\Orb(Y)\in\Ocal_d^1$. It is then very easy to see that the set of orbits of $\Xbf_n$ whose elements have the stabilizer
$\langle\sigma^d\rangle$ is
$$
\bigg\{\bigcup_{Y\in O} Y^{(n/d)}:O\in\Ocal_d^1\bigg\},
$$
which shows that the number of orbits of $\Xbf_n$ whose elements have the stabilizer $\langle\sigma^d\rangle$ is given by $\orb_1(d)$. Hence the
result.
\end{proof}

The next proposition shows that determining the sequences $\orb_d(n)$ reduces to determining only the two sequences $\orb_1(n)$ and $\orb_2(n)$.

\begin{proposition}\label{prop:Red}
For any $d|2n$, we have
$$
\orb_d(n) =
\begin{cases}
\orb_1(n/d), & \mbox{if $d$ is odd},\\
\orb_2(2n/d), & \mbox{if $d$ is even.}
\end{cases}
$$
\end{proposition}

\begin{proof}
Let $d\geqslant 1$ be an odd integer. Consider $O\in\Ocal_{dn}^d$ and $Y\in O$. Then
$$
\Stab(Y)=\langle\sigma^n\rangle=\{1,\sigma^n,\sigma^{2n},\ldots,\sigma^{(d-1)n}\}.
$$
By proceeding as in the proof of Proposition~\ref{prop:NOrbIncl}, we can easily see that
$$
\Ocal_{dn}^d =\bigg\{\bigcup_{X\in O} X^{(d)}:O\in\Ocal_n^1\bigg\},
$$
which obviously entails $\orb_d(dn)=\orb_1(n)$.

Now, let $d\geqslant 1$ be an even integer. Consider $O\in\Ocal_{dn/2}^d$ and $Y\in O$. Then either
$$
\Stab(Y)=\langle\sigma^{\frac n2}\rangle=\{1,\sigma^{\frac n2},\ldots,\sigma^{(d-1)\frac n2}\}
$$
(assuming $n$ even) or there is $j\in\{0,1,\ldots,dn/2-1\}$ such that
\begin{eqnarray*}
\Stab(Y) &=& \langle\sigma^n,\sigma^j\tau\rangle\\
&=& \{1,\sigma^n,\ldots,\sigma^{(\frac d2-1)n},\sigma^j\tau,\sigma^{n+j}\tau,\ldots,\sigma^{(\frac d2-1)n+j}\tau\}.
\end{eqnarray*}
In both cases, there is $X\in\Xbf_n$ such that $Y=X^{(d/2)}$ and, by Proposition~\ref{prop:XXd}, we have
$$
\Stab(X)=\{1,\sigma^{\frac n2}\},
$$
in the first case, and
$$
\Stab(X)=\{1,\sigma^{j\,(\mathrm{mod}~n)}\tau\},
$$
in the second case. It is then very easy to see that
$$
\Ocal_{dn/2}^d =\bigg\{\bigcup_{X\in O} X^{(d/2)}:O\in\Ocal_n^2\bigg\},
$$
which obviously entails $\orb_d(dn/2)=\orb_2(n)$.
\end{proof}

\section{Counting non-isomorphic MISs of $C_n$}\label{sec:MISCn}

In the present section we derive formulas for the enumeration of non-isomorphic MISs of $C_n$. Thus, we now assume that $\Xbf_n$ is the set of
MISs of $C_n$.

We can readily see that the membership function $\mathbf{1}_X$ of any MIS $X$ of $C_n$ is characterized by the following two conditions:
\begin{eqnarray}\label{eq:minmax}
\min\{\mathbf{1}_X(v_i),\mathbf{1}_X(v_{i+1})\}                     &=& 0 \qquad (i\in\Z),\label{eq:min}\\
\max\{\mathbf{1}_X(v_i),\mathbf{1}_X(v_{i+1}),\mathbf{1}_X(v_{i+2})\} &=& 1 \qquad (i\in\Z).\label{eq:max}
\end{eqnarray}
In fact, condition (\ref{eq:min}) expresses independence whereas condition (\ref{eq:max}) expresses maximality. Thus, for any MIS $X$ of $C_n$,
the $n$-tuple $\mathbf{1}_X(V_n)$ is a cyclic $n$-list made up of 0s and 1s, with one or two 0s between two neighboring 1s.

As mentioned in the introduction, we have $|\Xbf_n|=p(n)$, the $n$th term of the Perrin sequence. Then, from (\ref{eq:XDecomposition}), we
immediately have
\begin{equation}\label{eq:PerrinDecomposition}
p(n) = \sum_{d|2n} \frac{2n}d\,\orb_d(n).
\end{equation}

\begin{remark}
The fact that the Perrin sequence counts the distinct MISs of $C_n$ was observed without proof by F\"uredi \cite{Fur87}. A very simple bijective
proof, suggested by Vatter \cite{Vat}, can be written as follows. Let $X\in\Xbf_n$ and suppose that the vertex in $X$ of maximal label is $k$,
and the vertex of second-greatest label is $j$. Then $k$ must be $j+2$ or $j+3$. If $k=j+2$, then $X\setminus\{v_k\}$ can be viewed as a MIS of
$C_{n-2}$, while if $k=j+3$ then $X\setminus\{v_k\}$ can be viewed as a MIS of $C_{n-3}$. The inverse to this map can easily be described, thus
proving that the number of MISs of $C_n$ satisfies the Perrin recurrence.\qed
\end{remark}

We can readily see that any orbit of $\mathbf{X}_n$ can be uniquely represented by a cyclic list made up of $2$s and/or $3$s summing up to $n$,
where a clockwise writing is not distinguished from its counterclockwise counterpart. The bijection between this representation and the $n$-list
representation is straightforward: put a $2$ whenever only one $0$ separates two neighboring $1$s and put a $3$ whenever two $0$s separate two
neighboring $1$s.

This representation of orbits immediately leads to the following result:

\begin{proposition}\label{prop:OrbCy}
$\orb(n)$ is the number of cyclic compositions of $n$ in which each term is either $2$ or $3$, where a clockwise writing is not distinguished
from its counterclockwise counterpart.
\end{proposition}

Consider the Padovan sequence $q=(q(n))_{n\in\N}$ (shifted Sloane's A000931 \cite{Slo}), which is defined as $q(1)=0$, $q(2)=1$, $q(3)=1$, and
\begin{equation}\label{eq:PadovRec}
q(n)=q(n-2)+q(n-3)\qquad (n\geqslant 4),
\end{equation}
and let $r=(r(n))_{n\in\N}$ be the sequence defined by
$$
r(n)=
\begin{cases}
q(k), & \mbox{if $n=2k-1$},\\%
q(k+2), & \mbox{if $n=2k$}.
\end{cases}
$$
Note that, from (\ref{eq:PadovRec}) it follows that $r(n)$ fulfills the recurrence equation
$$
r(n)=r(n-4)+r(n-6)\qquad (n\geqslant 7).
$$

The following immediate proposition shows that $q(n)$ is the number of ways of writing $n$ as an ordered sum in which each term is either $2$ or
$3$. For example, $q(8) = 4$, and there are $4$ ways of writing $8$ as an ordered sum of $2$s and $3$s: $2 + 2 + 2 + 2$, $2 + 3 + 3$, $3 + 2 +
3$, $3 + 3 + 2$.

\begin{proposition}\label{proposition:CompQ}
$q(n)$ is the number of compositions of $n$ in which each term is either $2$ or $3$.
\end{proposition}

\begin{proof}
The values of $q(n)$ for $n\leqslant 3$ can be easily calculated. Now, fix $n\geqslant 4$. From among all the possible compositions of $n$,
$q(n-2)$ of them begin with a 2 and $q(n-3)$ of them begin with a 3, which establishes (\ref{eq:PadovRec}).
\end{proof}

As the following proposition shows, $r(n)$ shares the same property as $q(n)$ but with the additional condition that the terms form a cyclic and
palindromic composition. For example, $r(8) = 2$, and there are only $2$ ways of writing $8$ as an ordered sum of $2$s and $3$s in a cyclic and
palindromic way: $2 + 2 + 2 + 2$, $3 + 2 + 3$.

\begin{proposition}\label{proposition:CyPal}
$r(n)$ is the number of cyclic and palindromic compositions of $n$ in which each term is either $2$ or $3$.
\end{proposition}

\begin{proof}
Let $T\in\{2,3\}^m$ be a cyclic and palindromic composition of $n$, with $m$ terms placed counterclockwise and uniformly on the unit circle. Let
$X$ be a representative of the corresponding orbit of $\mathbf{X}_n$. Clearly, each symmetry axis of $X$ induces a symmetry axis of the cyclic
list $T$. The number of symmetry axes, say $d\geqslant 1$, necessarily divides $n$ since $\Stab(X)$ is a dihedral subgroup of $D_{2n}$.

Each of the $d$ axes intersects the unit circle at two points. Label these $2d$ intersection points as $a_0,a_1,\ldots,a_{2d-1}$
counterclockwise. These points constitute the vertices of a regular $2d$-gon. Let $S$ be the sublist of $T$ consisting of the terms located on
the arc $[a_0,a_1]$ inclusive, with the property that the terms located at $a_0$ or $a_1$, if any, are divided by two, and let $\overline{S}$ be
the reversed version of $S$. According to the kaleidoscope principle (see for instance \cite{Goo04} for an expository note on dihedral
kaleidoscopes and Coxeter groups), the sublist $S$ is located on each of the $d$ arcs of the form $[a_{2i},a_{2i+1}]$ $(i\in\Z_d)$ whereas the
sublist $\overline{S}$ is located on each of the $d$ arcs of the form $[a_{2i+1},a_{2i+2}]$ $(i\in\Z_d)$. In particular, the elements of
$S\cup\overline{S}$ sum up to $n/d$.

Suppose first that $n=2k-1$ is odd, which necessarily implies that $d$ is odd (since $d|n$). By symmetry, there is a point $a_i$ ($i\in\Z_{2d}$)
at which we have a $3$ and, by rotating the composition, we can assume that $i=0$. At the opposite point $a_d$ we have either a $2$ or no term.
As this feature holds for every symmetry axis, we immediately observe that the only possible palindromic writing of $T$ whose two extreme terms
are $\frac 32$s is obtained by the alternating concatenation of $2d$ sublists $S$ and $\overline{S}$, that is,
$$
W = S\overline{S}S\overline{S}\cdots S\overline{S}.
$$
Let $L$ be the first half of $W$ and, except for its first element (which is $\frac 32$) and its last element (which can be $1$, $2$, or $3$),
replace in $L$ each sublist $(\frac 32,\frac 32)$ with $3$ and each sublist $(1,1)$ with $2$.

Let $\mathbf{T}_n$ be the set of cyclic and palindromic compositions of $n$ in which each term is either $2$ or $3$. Clearly, we can assume that
each of these compositions is written so that its two extreme terms are $\frac 32$s. For any $T\in\mathbf{T}_n$, denote by $L^T$ the sublist as
defined above. Then, it is easy to see that the disjoint union
$$
\bigcup_{T\in\mathbf{T}_n}L^T,
$$
which has cardinality $|\mathbf{T}_n|$, is also the set of compositions of $\frac 32+\frac{n-3}2=\frac 32+(k-2)$ in which each term is either
$2$ or $3$, except for the first term, which is $\frac 32$, and the last term which can be $1$, or $2$, or $3$. By calculating the cardinality
of this latter set, we obtain
$$
|\mathbf{T}_n| = q(k-2)+q(k-3).
$$
Indeed, by Proposition~\ref{proposition:CompQ}, there are $q(k-2)$ compositions whose last term is either $2$ or $3$, plus $q(k-3)$ compositions
whose last term is $1$. Finally, using (\ref{eq:PadovRec}) we obtain $|\mathbf{T}_n|=q(k)=r(n)$, which proves the result for $n=2k-1$.

Let us now assume that $n=2k$ is even. Then the two possible palindromic writings of $T$ are obtained by the alternating concatenation of $2d$
sublists $S$ and $\overline{S}$, that is,
\begin{eqnarray*}
W_1 &=& S\overline{S}S\overline{S}\cdots S\overline{S},\\
W_2 &=& \overline{S}S\overline{S}S\cdots \overline{S}S.
\end{eqnarray*}
(Note that if $S$ contains at least two elements then it cannot be palindromic for otherwise we would have more than $d$ symmetry axes.)

Let $L_1$ (resp.\ $L_2$) be the first half of $W_1$ (resp.\ $W_2$) and, without modifying its two extreme elements, replace in $L_1$ (resp.\
$L_2$) each sublist $(\frac 32,\frac 32)$ with $3$ and each sublist $(1,1)$ with $2$. If $S$ contains at least two elements then $L_1$ and $L_2$
are two distinct sublists. If $S$ contains only one element then we simply set $L_1=(1,2,\ldots,2,1)$ and $L_2=(2,\ldots,2)$ in case $S=(1)$ and
$L_1=(\frac 32,3,\ldots,3,\frac 32)$ and $L_2=(3,\ldots,3)$ in case $S=(\frac 32)$.

Let $\mathbf{T}_n$ be the set of cyclic and palindromic compositions of $n$ in which each term is either $2$ or $3$ and, for any
$T\in\mathbf{T}_n$, denote by $L_1^T$ and $L_2^T$ the two sublists as defined above. Then, it is easy to see that the disjoint union
$$
\bigcup_{T\in\mathbf{T}_n}L^T_1\cup\bigcup_{T\in\mathbf{T}_n}L^T_2,
$$
which has cardinality $2|\mathbf{T}_n|$, is also the set of compositions of $\frac n2=k$ in which each term is either $2$ or $3$, except for the
two extreme terms which can be $1$, or $\frac 32$, or $2$, or $3$. By calculating the cardinality of this latter set, we obtain
\begin{eqnarray*}
2\,|\mathbf{T}_n| &=& q(k)+2q(k-1)+q(k-2)+q(k-3)\\
&=& 2q(k)+2q(k-1)\\
&=& 2q(k+2).
\end{eqnarray*}
Indeed, by Proposition~\ref{proposition:CompQ}, there are $q(k)$ compositions whose two extreme terms are $2$s or $3$s, plus $2q(k-1)$
compositions whose only one extreme term is $1$, plus $q(k-2)$ compositions whose two extreme terms are $1$s, plus $q(k-3)$ compositions whose
two extreme terms are $\frac 32$s. Finally, we have $|\mathbf{T}_n|=q(k+2)=r(n)$, which completes the proof.
\end{proof}

The following example illustrates the proof of Proposition~\ref{proposition:CyPal} in the case $n=16$.

\begin{example}
There are $7$ cyclic and palindromic compositions of $16$ in which each term is either $2$ or $3$. Table~\ref{tab:ex} gives these compositions
together with the number $d$ of symmetry axes and the sublists $L_1$ and $L_2$, as defined in the proof of Proposition~\ref{proposition:CyPal}.
\begin{table}[ht]
$$
\begin{array}{|c|c|c|c|}
\hline T & d & L_1 & L_2 \\
\hline %
(2,2,2,2,2,2,2,2) & 8 & (2,2,2,2) & (1,2,2,2,1) \\
(3,2,3,3,2,3) & 2 & (3,2,3) & (1,3,3,1) \\
(2,2,3,2,3,2,2) & 1 & (2,2,3,1) & (1,3,2,2) \\
(2,3,2,2,2,3,2) & 1 & (2,3,2,1) & (1,2,3,2) \\
(2,3,3,3,3,2) & 1 & (2,3,3) & (3,3,2) \\
(3,2,2,2,2,2,3) & 1 & (3,2,2,1) & (1,2,2,3) \\
(\frac 32,2,3,3,3,2,\frac 32) & 1 & (\frac 32,2,3,\frac 32) & (\frac 32,3,2,\frac 32)\\
\hline
\end{array}
$$
\caption{Cyclic and palindromic compositions of $16$} \label{tab:ex}
\end{table}
\end{example}

Through the identification of the orbits of $\mathbf{X}_n$ whose elements have at least one symmetry axis with the cyclic and palindromic
compositions of $n$ made up of $2$s and/or $3$s, Proposition~\ref{proposition:CyPal} can be immediately rewritten as follows:

\begin{proposition}\label{prop:NOrbNIncl}
The number of orbits of $\Xbf_n$ whose elements have a stabilizer not included in $\langle\sigma\rangle$ is given by $r(n)$.
\end{proposition}

We are now able to state our main result, which gives explicit expressions for the sequences $\orb$, $\orb_1$, and $\orb_2$. To this end, recall
that the \emph{Dirichlet convolution product}\/ of two sequences $f=(f(n))_{n\in\N}$ and $g=(g(n))_{n\in\N}$ is the sequence $f\ast g=((f\ast
g)(n))_{n\in\N}$ defined as
$$
(f\ast g)(n)=\sum_{d|n}f(d)\, g(\frac nd).
$$
Moreover, in addition to the sequences $p(n)$, $q(n)$, and $r(n)$ introduced above, consider the following integer sequences:
\begin{eqnarray*}
\mu(n) &=& \mathrm{A008683}(n)\qquad \mbox{(\emph{M\"obius function})},\\
\mathrm{A113788}(n) &=& \frac 1n\, (p\ast\mu)(n),\\
\mathbf{1}(n) &=& 1,\\
e_k(n) &=& \begin{cases} 1, & \mbox{if $n=k$,}\\ 0, & \mbox{else,}\end{cases}\qquad (k\in\N).
\end{eqnarray*}

\begin{theorem}\label{thm:main}
There holds
\begin{eqnarray}
\orb &=& r+\orb_1\ast\mathbf{1},\label{eq:orb}\\
2\orb_1 &=& \mathrm{A113788}-r\ast\mu,\label{eq:orb1}\\
\orb_2 &=& r\ast\mu + \orb_1\ast e_2.\label{eq:orb2}
\end{eqnarray}
\end{theorem}

\begin{proof}
Combining Propositions~\ref{prop:NOrbIncl} and \ref{prop:NOrbNIncl} immediately leads to (\ref{eq:orb}). Next, combining (\ref{eq:NumberOrbits})
and (\ref{eq:orb}), we can write
$$
r(n) = \sum_{d|2n}\orb_d(n)-\sum_{d|n}\orb_1(d),
$$
that is, using Proposition~\ref{prop:Red},
\begin{eqnarray*}
r(n)
&=& \sum_{\textstyle{d|2n\atop d~\mathrm{even}}}\orb_d(n)+\sum_{\textstyle{d|2n\atop d~\mathrm{odd}}}\orb_d(n)-\sum_{d|n}\orb_1(d)\\
&=& \sum_{\textstyle{d|2n\atop \frac{2n}d~\mathrm{even}}}\orb_2(d)+\sum_{\textstyle{d|n\atop \frac nd~\mathrm{odd}}}\orb_1(d)-\sum_{d|n}\orb_1(d)\\
&=& \sum_{d|n}\orb_2(d)-\sum_{\textstyle{d|n\atop \frac nd~\mathrm{even}}}\orb_1(d).
\end{eqnarray*}
Hence, introducing $g_2(n)=\frac{1+(-1)^n}2$, we obtain
$$
r=\orb_2\ast \mathbf{1}-\orb_1\ast g_2,
$$
that is, using M\"obius inversion formula and the immediate identity $g_2=e_2\ast\mathbf{1}$,
\begin{equation}\label{eq:ro201}
r\ast\mu=\orb_2-\orb_1\ast e_2.
\end{equation}

On the other hand, combining (\ref{eq:PerrinDecomposition}) with Proposition~\ref{prop:Red}, we obtain
\begin{eqnarray*}
p(n) &=& \sum_{\textstyle{d|2n\atop d~\mathrm{odd}}} \frac{2n}d \, \orb_1\big(\frac nd\big)
+\sum_{\textstyle{d|2n\atop d~\mathrm{even}}} \frac{2n}d \,\orb_2\big(\frac{2n}d\big)\\
&=& \sum_{\textstyle{d|n\atop d~\mathrm{odd}}} \frac{2n}d \, \orb_1\big(\frac nd\big)
+\sum_{\textstyle{d|2n\atop \frac{2n}d~\mathrm{even}}} d \,\orb_2(d)\\
&=& \sum_{\textstyle{d|n\atop \frac nd~\mathrm{odd}}} 2d \, \orb_1(d) +\sum_{d|n} d \,\orb_2(d).
\end{eqnarray*}
Hence, introducing $f_1(n)=n\,\orb_1(n)$, $f_2(n)=n\,\orb_2(n)$, and $g_1(n)=\frac{1-(-1)^n}2$, we obtain
$$
p=2f_1\ast g_1+f_2\ast \mathbf{1},
$$
that is, using M\"obius inversion formula and the fact that $g_1=\mathbf{1}-g_2=\mathbf{1}-e_2\ast\mathbf{1}$,
$$
p\ast\mu=2f_1-2f_1\ast e_2+f_2,
$$
which implies
\begin{equation}\label{eq:po201}
\mathrm{A113788} = 2\orb_1-\orb_1\ast e_2+\orb_2.
\end{equation}
Combining (\ref{eq:ro201}) and (\ref{eq:po201}) immediately leads to formulas (\ref{eq:orb1}) and (\ref{eq:orb2}).
\end{proof}

It is noteworthy that the sequence $\orb(n)$ can be rewritten as
$$
2\orb(n)=r(n)+\frac 1n(p\ast\phi)(n),
$$
where $\phi=\mathrm{A000010}$ is the \emph{Euler totient function}. In fact, combining (\ref{eq:orb}) and (\ref{eq:orb1}), we obtain
\begin{eqnarray*}
2\orb(n) &=& 2r(n)+2(\orb_1\ast\mathbf{1})(n)\\
&=& r(n)+(\mathrm{A113788}\ast\mathbf{1})(n)
\end{eqnarray*}
and the latter term also writes
\begin{equation}\label{eq:PastPhi}
(\mathrm{A113788}\ast\mathbf{1})(n) = \sum_{d|n} \frac dn(p\ast\mu)\big(\frac nd\big) = \frac 1n(p\ast\mu\ast\mathrm{Id})(n)=\frac
1n(p\ast\phi)(n),
\end{equation}
where $\mathrm{Id}(n)=n$ is the identity sequence.

\section{Counting unlabeled MISs of $C_n$}

We now investigate the simpler problem of enumerating the \emph{unlabeled}\/ MISs of $C_n$, that is, the MISs defined up to a rotation. In fact,
these unlabeled MISs are nothing less than the orbits of $\mathbf{X}_n$ under the action of $\langle\sigma\rangle$.

Let $\orb^{\sigma}(n)$ denote the number of unlabeled MISs of $C_n$ and, for any $d|n$, let $\orb_d^{\sigma}(n)$ denote the number of unlabeled
MISs of $C_n$ that have $\frac nd$ isomorphic representatives or, equivalently, those that have a stabilizer of cardinality $d$.

Then, we clearly have
\begin{eqnarray}
\orb^{\sigma}(n) &=& \sum_{d|n} \orb_d^{\sigma}(n),\label{eq:NumberOrbitsSigma}\\
p(n) &=& \sum_{d|n} \frac nd\,\orb_d^{\sigma}(n),\label{eq:PerrinDecompositionSigma}
\end{eqnarray}
and, by proceeding as in the proof of Proposition~\ref{prop:Red}, we can easily show that
\begin{equation}\label{eq:RedSigma}
\orb_d^{\sigma}(n) = \orb_1^{\sigma}(n/d)\qquad (d|n).
\end{equation}

We then obtain the following formulas:

\begin{theorem}\label{thm:mainSigma}
There holds
\begin{eqnarray}
\orb^{\sigma} &=& \orb_1^{\sigma}\ast\mathbf{1},\label{eq:orbSigma}\\
\orb_1^{\sigma} &=& \mathrm{A113788}.\label{eq:orb1Sigma}
\end{eqnarray}
\end{theorem}

\begin{proof}
Define the sequence $f_1^{\sigma}(n)=n\,\orb_1^{\sigma}(n)$. Combining (\ref{eq:NumberOrbitsSigma}) and (\ref{eq:RedSigma}) immediately leads to
(\ref{eq:orbSigma}). On the other hand, combining (\ref{eq:PerrinDecompositionSigma}) and (\ref{eq:RedSigma}), we obtain
$$
p(n)=\sum_{d|n} \frac nd\,\orb_1^{\sigma}(n/d)=(f_1^{\sigma}\ast\mathbf{1})(n)
$$
and by using M\"obius inversion formula we are immediately led to (\ref{eq:orb1Sigma}).
\end{proof}

It is noteworthy that, combining (\ref{eq:PastPhi}), (\ref{eq:orbSigma}), and (\ref{eq:orb1Sigma}), we also obtain
$$
\orb^{\sigma}(n)=\frac 1n(p\ast\phi)(n).
$$

Finally, we also have the immediate result:

\begin{proposition}\label{prop:OrbSCy}
$\orb^{\sigma}(n)$ is the number of cyclic compositions of $n$ in which each term is either $2$ or $3$.
\end{proposition}

\begin{remark}
Just as for the results obtained in Section~\ref{sec:autCn}, the results stated in the present section do not essentially lie on the intrinsic
properties of MISs. Up to the left-hand side of (\ref{eq:PerrinDecompositionSigma}), formulas (\ref{eq:NumberOrbitsSigma})--(\ref{eq:RedSigma})
remain valid for any other definition of $\mathbf{X}_n$.
\end{remark}

\section{Summary}

Starting from Perrin (A001608) and Padovan (shifted A000931) sequences, respectively denoted $p(n)$ and $q(n)$ in the present paper, we have
introduced the following new integer sequences, with an explicit expression for each of them:

\begin{itemize}
\item $r(n)=\mathrm{A127682}(n)$ is the number of non-isomorphic MISs of $C_n$ having at least one symmetry axis (see
Proposition~\ref{prop:NOrbNIncl}), where two MISs are isomorphic if they are identical up to a rotation and a reflection. It is also the number
of cyclic and palindromic compositions of $n$ in which each term is either $2$ or $3$ (see Proposition~\ref{proposition:CyPal}). Recall that
this sequence is defined as
$$
r(n)=
\begin{cases}
q(k), & \mbox{if $n=2k-1$},\\%
q(k+2), & \mbox{if $n=2k$}.
\end{cases}
$$

\item For any $d|2n$, $\orb_d(n)$ gives the number of non-isomorphic MISs of $C_n$ having $\frac{2n}d$ isomorphic representatives. This sequence
can always be expressed from one of the sequences $\orb_1(n)=\mathrm{A127683}(n)$ and $\orb_2(n)=\mathrm{A127684}(n)$ (see
Proposition~\ref{prop:Red}), which in turn can be directly calculated from the formulas
\begin{eqnarray*}
\orb_1(n) &=& \frac 12\big(\mathrm{A113788}(n)-(r\ast\mu)(n)\big)\\
\orb_2(n) &=& (r\ast\mu)(n)+(\orb_1\ast e_2)(n).
\end{eqnarray*}

\item $\orb(n)=\mathrm{A127685}(n)$ gives the number of non-isomorphic MISs of $C_n$. It is also the number of cyclic compositions of $n$ in
which each term is either $2$ or $3$, where a clockwise writing is not distinguished from its counterclockwise counterpart (see
Proposition~\ref{prop:OrbCy}). This sequence is given explicitly by
$$
\orb(n)=r(n)+(\orb_1\ast\mathbf{1})(n)=\frac 12\Big(r(n)+\frac 1n(p\ast\phi)(n)\Big).
$$

\item $(\orb_1\ast\mathbf{1})(n)=\mathrm{A127686}(n)$ is the number of non-isomorphic MISs of $C_n$ having no symmetry axis (see
Proposition~\ref{prop:NOrbIncl}). It is also the number of cyclic and non-palindromic compositions of $n$ in which each term is either $2$ or
$3$, where a clockwise writing is not distinguished from its counterclockwise counterpart. This sequence is given explicitly by
$$
(\orb_1\ast\mathbf{1})(n)=\orb(n)-r(n).
$$

\item $\orb^{\sigma}(n)=\mathrm{A127687}(n)$ is the number of unlabeled MISs of $C_n$. It is also the number of cyclic compositions of $n$ in
which each term is either $2$ or $3$ (see Proposition~\ref{prop:OrbSCy}). This sequence is given explicitly by
$$
\orb^{\sigma}(n)=(\mathrm{A113788}\ast\mathbf{1})(n) = \frac 1n(p\ast\phi)(n) = 2\orb(n)-r(n).
$$

\item For any $d|n$, $\orb_d^{\sigma}(n)$ gives the number of unlabeled MISs of $C_n$ having $\frac nd$ isomorphic representatives. This
sequence can always be expressed from the sequence $\orb_1^{\sigma}(n)$ (see (\ref{eq:RedSigma})), which in turn can be directly calculated from
the formulas
$$
\orb_1^{\sigma}(n)=\mathrm{A113788}(n)=(\orb^{\sigma}\ast\mu)(n).
$$
\end{itemize}

The first 40 values of the main sequences considered in this paper are listed in Table~\ref{tab:seq}. We observe that, when $n$ is prime, we
have $\orb(n)=\orb_1(n)+\orb_2(n)$, $\orb_2(n)=r(n)$, and $\orb^{\sigma}(n)=\orb_1^{\sigma}(n)=\orb(n)+\orb_1(n)=\frac 1n p(n)$, which can be
easily verified.

\section{A musical application}

In classical tonal music (see for instance Benward and Saker~\cite{BenSak08}), MISs on $C_{7k}$ for $k=1,2,\ldots$ may be seen as musical chords
built on thirds and fourths, where $k$ denotes the number of octaves on some major or minor scale. For instance, in C major the Cmaj chord,
namely C-E-G, is a MIS on $C_7$ and the thirteenth chord C13maj, namely C-E-G-B-d-f-a, is a MIS on $C_{14}$.

On $C_7$ there exists only $\orb(7)=\orb^{\sigma}(7)=1$ non-isomorphic MIS, with $p(7)=7$ elements giving all possible chords built on generic
thirds and fourths. However, on $C_{14}$, that is a two octave scale, we observe $\orb(14)=\orb^{\sigma}(14)=5$ non-isomorphic MIS, with
$p(14)=51$ elements giving all possible chords on generic (unlabeled) thirds and fourths.

\begin{figure}[htb]
\begin{center}
\includegraphics[height=.70\linewidth, keepaspectratio=true]{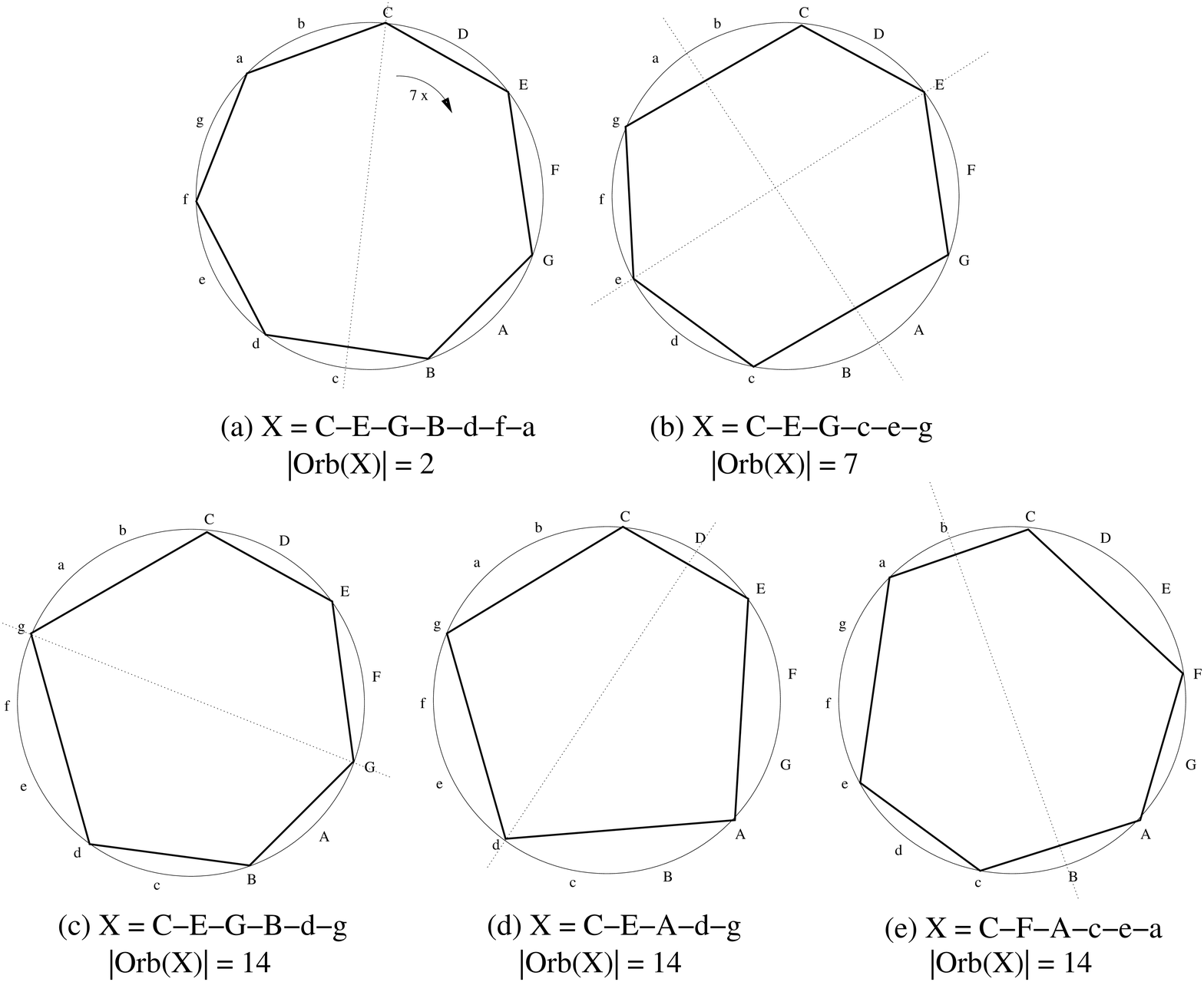}
\caption{Chord sets in C major built from generic thirds and fourths on two octaves} \label{fig:Chord}
\end{center}
\end{figure}

Figure~\ref{fig:Chord} shows the $5$ possible chord sets in C major built from generic thirds and fourths on two octaves. MIS (a) represents the
above-mentioned thirteenth chord C13maj, in which we observe $7$ symmetry axes passing through every note of the chord. In this case, the orbit
contains only $2$ representatives. MIS (b) represents a repetition of the Cmaj chord in each octave, in which we have $2$ symmetry axes, and the
orbit contains $7$ representatives. In the three remaining MISs (c), (d), and (e), we have only one symmetry axis passing respectively through
$2$, $1$, and $0$ notes and each of the corresponding orbits contains $14$ representatives.

Our results provide useful tools for analyzing the symmetries of third and fourth chords on more than two octaves.

\begin{table}[p]
$$
\begin{array}{|r|r|r|r|r|r|r|r|r|}
\hline %
n & p(n) & q(n) & r(n) & \orb(n) & \orb_1(n) & \orb_2(n) & \orb^{\sigma}(n) & \orb_1^{\sigma}(n) \\
\hline\hline %
1 & 0 & 0 & 0 & 0 & 0 & 0 & 0 & 0\\
2 & 2 & 1 & 1 & 1 & 0 & 1 & 1 & 1\\
3 & 3 & 1 & 1 & 1 & 0 & 1 & 1 & 1\\
4 & 2 & 1 & 1 & 1 & 0 & 0 & 1 & 0\\
5 & 5 & 2 & 1 & 1 & 0 & 1 & 1 & 1\\
6 & 5 & 2 & 2 & 2 & 0 & 0 & 2 & 0\\
7 & 7 & 3 & 1 & 1 & 0 & 1 & 1 & 1\\
8 & 10 & 4 & 2 & 2 & 0 & 1 & 2 & 1\\
9 & 12 & 5 & 2 & 2 & 0 & 1 & 2 & 1\\
10 & 17 & 7 & 3 & 3 & 0 & 1 & 3 & 1\\
\hline 11 & 22 & 9 & 2 & 2 & 0 & 2 & 2 & 2\\
12 & 29 & 12 & 4 & 4 & 0 & 2 & 4 & 2\\
13 & 39 & 16 & 3 & 3 & 0 & 3 & 3 & 3\\
14 & 51 & 21 & 5 & 5 & 0 & 3 & 5 & 3\\
15 & 68 & 28 & 4 & 5 & 1 & 2 & 6 & 4\\
16 & 90 & 37 & 7 & 7 & 0 & 5 & 7 & 5\\
17 & 119 & 49 & 5 & 6 & 1 & 5 & 7 & 7\\
18 & 158 & 65 & 9 & 10 & 1 & 6 & 11 & 8\\
19 & 209 & 86 & 7 & 9 & 2 & 7 & 11 & 11\\
20 & 277 & 114 & 12 & 14 & 2 & 9 & 16 & 13\\
\hline 21 & 367 & 151 & 9 & 14 & 5 & 7 & 19 & 17\\
22 & 486 & 200 & 16 & 20 & 4 & 13 & 24 & 21\\
23 & 644 & 265 & 12 & 20 & 8 & 12 & 28 & 28\\
24 & 853 & 351 & 21 & 30 & 9 & 16 & 39 & 34\\
25 & 1130 & 465 & 16 & 31 & 15 & 15 & 46 & 45\\
26 & 1497 & 616 & 28 & 44 & 16 & 24 & 60 & 56\\
27 & 1983 & 816 & 21 & 48 & 27 & 19 & 75 & 73\\
28 & 2627 & 1081 & 37 & 67 & 30 & 32 & 97 & 92\\
29 & 3480 & 1432 & 28 & 74 & 46 & 28 & 120 & 120\\
30 & 4610 & 1897 & 49 & 104 & 54 & 44 & 159 & 151\\
\hline 31 & 6107 & 2513 & 37 & 117 & 80 & 37 & 197 & 197\\
32 & 8090 & 3329 & 65 & 161 & 96 & 58 & 257 & 250\\
33 & 10717 & 4410 & 49 & 188 & 139 & 46 & 327 & 324\\
34 & 14197 & 5842 & 86 & 254 & 167 & 81 & 422 & 414\\
35 & 18807 & 7739 & 65 & 302 & 237 & 63 & 539 & 537\\
36 & 24914 & 10252 & 114 & 407 & 292 & 104 & 700 & 687\\
37 & 33004 & 13581 & 86 & 489 & 403 & 86 & 892 & 892\\
38 & 43721 & 17991 & 151 & 654 & 501 & 145 & 1157 & 1145\\
39 & 57918 & 23833 & 114 & 801 & 687 & 110 & 1488 & 1484\\
40 & 76725 & 31572 & 200 & 1064 & 862 & 189 & 1928 & 1911\\
\hline
\end{array}
$$
\caption{First 40 values of the main sequences} \label{tab:seq}
\end{table}

\end{document}